\documentclass[a4paper]{article}
\usepackage{amsmath}
\usepackage{amssymb}
\font\black=cmbx10 \font\sblack=cmbx7 \font\ssblack=cmbx5
\font\blackital=cmmib10  \skewchar\blackital='177 \font\sblackital=cmmib7
\skewchar\sblackital='177 \font\ssblackital=cmmib5
\skewchar\ssblackital='177 \font\sanss=cmss10 \font\ssanss=cmss8 scaled
900 \font\sssanss=cmss8 scaled 600 \font\blackboard=msbm10
\font\sblackboard=msbm7 \font\ssblackboard=msbm5 \font\caligr=eusm10
\font\scaligr=eusm7 \font\sscaligr=eusm5 
\font\fraktur=eufm10 \font\sfraktur=eufm7 \font\ssfraktur=eufm5

\font\bsymb=cmsy10 scaled\magstep2
\def\all#1{\setbox0=\hbox{\lower1.5pt\hbox{\bsymb
       \char"38}}\setbox1=\hbox{$_{#1}$} \box0\lower2pt\box1\;}
\def\exi#1{\setbox0=\hbox{\lower1.5pt\hbox{\bsymb \char"39}}
       \setbox1=\hbox{$_{#1}$} \box0\lower2pt\box1\;}

\def\tx#1{{\fam0\relax#1}}

\newfam\bifam
\textfont\bifam=\blackital \scriptfont\bifam=\sblackital
\scriptscriptfont\bifam=\ssblackital

\newfam\blfam
\textfont\blfam=\black \scriptfont\blfam=\sblack
\scriptscriptfont\blfam=\ssblack

\newfam\bbfam
\textfont\bbfam=\blackboard \scriptfont\bbfam=\sblackboard
\scriptscriptfont\bbfam=\ssblackboard

\newfam\ssfam
\textfont\ssfam=\sanss \scriptfont\ssfam=\ssanss
\scriptscriptfont\ssfam=\sssanss
\def\sss#1{{\fam\ssfam\relax#1}}

\newfam\clfam
\textfont\clfam=\caligr \scriptfont\clfam=\scaligr
\scriptscriptfont\clfam=\sscaligr

\newfam\frfam
\textfont\frfam=\fraktur \scriptfont\frfam=\sfraktur
\scriptscriptfont\frfam=\ssfraktur

\def\hpb#1{\setbox0=\hbox{${#1}$}
    \copy0 \kern-\wd0 \kern.2pt \box0}
\def\vpb#1{\setbox0=\hbox{${#1}$}
    \copy0 \kern-\wd0 \raise.08pt \box0}

\def\pmb#1{\setbox0\hbox{${#1}$} \copy0 \kern-\wd0 \kern.2pt \box0}
\def\pmbb#1{\setbox0\hbox{${#1}$} \copy0 \kern-\wd0
      \kern.2pt \copy0 \kern-\wd0 \kern.2pt \box0}
\def\pmbbb#1{\setbox0\hbox{${#1}$} \copy0 \kern-\wd0
      \kern.2pt \copy0 \kern-\wd0 \kern.2pt
    \copy0 \kern-\wd0 \kern.2pt \box0}
\def\pmxb#1{\setbox0\hbox{${#1}$} \copy0 \kern-\wd0
      \kern.2pt \copy0 \kern-\wd0 \kern.2pt
      \copy0 \kern-\wd0 \kern.2pt \copy0 \kern-\wd0 \kern.2pt \box0}
\def\pmxbb#1{\setbox0\hbox{${#1}$} \copy0 \kern-\wd0 \kern.2pt
      \copy0 \kern-\wd0 \kern.2pt
      \copy0 \kern-\wd0 \kern.2pt \copy0 \kern-\wd0 \kern.2pt
      \copy0 \kern-\wd0 \kern.2pt \box0}

\def\sT{{\sss T}}

\def\xi{\tx{i}}

\newcommand{\be}{\begin{equation}}
\newcommand{\ee}{\end{equation}}
\newcommand{\ra}{\rightarrow}

\newcommand{\bea}{\begin{eqnarray}}
\newcommand{\eea}{\end{eqnarray}}
\newcommand{\beas}{\begin{eqnarray*}}
\newcommand{\eeas}{\end{eqnarray*}}

\newcommand{\R}{\mathbb{R}}

\newcommand{\nn}{\nonumber}
\newcommand{\ot}{\otimes}
\newcommand{\op}{\oplus}
\newcommand{\s}{\circ}

\newcommand{\pa}{\partial}

\newcommand{\A}{{\cal A}}
\newcommand{\Li}{{\cal L}}

\newcommand{\Ll}{\Li}

\def\la{\langle}
\def\ran{\rangle}

\def\Ci{C^\infty}

\mathchardef\za="710B  
\mathchardef\zb="710C  
\mathchardef\zg="710D  
\mathchardef\zd="710E  
\mathchardef\zve="710F 
\mathchardef\zz="7110  
\mathchardef\zh="7111  
\mathchardef\zvy="7112 
\mathchardef\zi="7113  
\mathchardef\zk="7114  
\mathchardef\zl="7115  
\mathchardef\zm="7116  
\mathchardef\zn="7117  
\mathchardef\zx="7118  
\mathchardef\zp="7119  
\mathchardef\zr="711A  
\mathchardef\zs="711B  
\mathchardef\zt="711C  
\mathchardef\zu="711D  
\mathchardef\zvf="711E 
\mathchardef\zq="711F  
\mathchardef\zc="7120  
\mathchardef\zw="7121  
\mathchardef\ze="7122  
\mathchardef\zy="7123  
\mathchardef\zf="7124  
\mathchardef\zvr="7125 
\mathchardef\zvs="7126 
\mathchardef\zf="7127  
\mathchardef\zG="7000  
\mathchardef\zD="7001  
\mathchardef\zY="7002  
\mathchardef\zL="7003  
\mathchardef\zX="7004  
\mathchardef\zP="7005  
\mathchardef\zS="7006  
\mathchardef\zU="7007  
\mathchardef\zF="7008  
\mathchardef\zW="700A  

\newcommand{\epf}{\hfill$\Box$}
\newcommand{\bepf}{\textit{Proof.-} }

\begin{document}

\title{Courant algebroid and Lie bialgebroid contractions}
\author{Jos\'e F. Cari\~nena\\
Depto. F\'{\i}sica Te\'orica, Universidad de Zaragoza\\
50009 Zaragoza, Spain\\
{\it e-mail:} jfc@posta.unizar.es
\and
Janusz Grabowski\thanks{Supported by KBN, grant No 2 P03A 020 24.}\\
Mathematical Institute, Polish Academy of Sciences\\ ul. \'Sniadeckich
8, P. O. Box 21, 00-956 Warszawa, Poland\\ {\it e-mail:}
jagrab@impan.gov.pl
\and
Giuseppe Marmo\\
Dipartimento di Scienze Fisiche,
Universit\`a Federico II di Napoli\\
and\\
INFN, Sezione di Napoli\\
Complesso Universitario di Monte Sant'Angelo\\
Via Cintia, 80125 Napoli, Italy\\
{\it e-mail:} marmo@na.infn.it}
\maketitle

\newtheorem{re}{Remark}
\newtheorem{theo}{Theorem}
\newtheorem{prop}{Proposition}
\newtheorem{lem}{Lemma}
\newtheorem{cor}{Corollary}
\newtheorem{ex}{Example}

\begin{abstract}
Contractions of Leibniz algebras and Courant algebroids by means of
(1,1)-tensors are introduced and studied. An appropriate version of
Nijenhuis tensors leads to natural deformations of Dirac structures and
Lie bialgebroids. One recovers presymplectic-Nijenhuis structures,
Poisson-Nijenhuis structures, and triangular Lie bialgebroids as
particular examples.

\bigskip\noindent
\textit{\textbf{MSC 2000:} Primary 17B99; Secondary 17B62, 53C15, 53D17.}

\medskip\noindent
\textit{\textbf{Key words:} Nijenhuis tensor, Lie bialgebroid, Leibniz
algebra, Courant bracket, Courant algebroid.}

\end{abstract}

\section{Introduction}
This note is a natural continuation of our previous work \cite{CGM}, where
contractions and Nijenhuis tensors have been studied for algebraic
operations of arbitrary type on sections of vector bundles. Recall that a
\textit{Nijenhuis tensor} $N$ for a bilinear operation "$\s$" on sections
of a vector bundle $A$ over $M$ is a $(1,1)$-tensor $N\in Sec(A\ot A^*)$
viewed as vector bundle morphism $N:A\ra A$ (or the corresponding
$\Ci(M)$-linear map $N:Sec(A)\ra Sec(A)$ on sections) such that its
\textit{Nijenhuis torsion}
\begin{equation}\label{1}
T_N(X,Y)=N(X)\s N(Y)-N(X\s_NY)
\end{equation}
vanishes, where "$\s_N$" is the contracted product:
\begin{equation}\label{2}
X\s_NY=N(X)\s Y+X\s N(Y)-N(X\s Y).
\end{equation}
The theory of Nijenhuis tensors for Lie algebra brackets goes back to a
concept of contractions of Lie algebras introduced by E.~J.~Saletan
\cite{Sa}. Nijenhuis tensors for Lie algebroids and Nijenhuis tensors on
Poisson manifolds were studied in \cite{MM,KSM} and in a number of
following papers. In \cite{CGM0} the authors of this note developed the
theory of Nijenhuis tensors for associative products, and in \cite{CGM} --
for arbitrary algebraic operations.

One can apply directly the procedures from \cite{CGM} to Leibniz algebras.
The vanishing of the Nijenhuis torsion implies that the contracted product
"$\s_N$" is again a Leibniz product. However, as we will see in the
example of the Courant product on $\sT M\op\sT^*M$ (in its Leibniz
version) the vanishing of the Nijenhuis torsion is a too restrictive
assumption. To get that "$\s_N$" is Leibniz it is sufficient to require
that $T_N$ is a Leibniz 2-cocycle. We will refer to such tensors $N$ as to
\textit{weak Nijenhuis tensors} for Leibniz algebras. Since the use of
weak Nijenhuis tensors does not lead to contractions in the strict sense
(they do not come from a limit procedure), one should rather call "$\s_N$"
a \textit{deformed product} in this case. So the convention throughout
this paper is that we use the word `contraction' heuristically, thinking
just on a procedure of passing from a product "$\s$" to the product
"$\s_N$" for a specifically chosen $(1,1)$-tensor $N$.

To introduce a notion of a Lie bialgebroid contraction we use the concept
of \textit{Courant algebroid} \cite{LWX} in its Leibniz version. Since the
Courant algebroid is not only a Leibniz product but also a non-degenerate
pairing with certain consistency conditions with the Leibniz product, we
check what property of $N$ ensures the consistency conditions being
satisfied also for "$\s_N$". It turns out that it is sufficient to assume
that $N+N^*=\zl I$, $\zl\in\R$, where $N^*$ is dual to $N$ with respect to
the pairing; we will call such tensors \textit{paired}. Thus, paired and
weak Niejnhuis tensors on Courant algebroids give rise to deformed Courant
algebroids.

There is a straightforward but very useful generalization of the concept
of the Nijenhuis tensor. Suppose that $L$ is a subbundle of $A$ whose
sections are closed with respect to the operation "$\s$", i.e. they form a
subalgebra in $(Sec(A),\s)$. If $Sec(L)$ is closed also for "$\s_N$" and
the torsion $T_N$ vanishes on $L$, i.e. $T_N(X,Y)=0$ for all $X,Y\in
Sec(L)$, we will refer to $N$ as to an \textit{outer Nijenhuis tensor} for
$(L,\s)$. This concept seems to be the right tool in contracting
\textit{Dirac structures}, i.e. subbundles of Courant algebroids which are
maximal isotropic and closed with respect to the product. In this approach
a \textit{Dirac-Nijenhuis} structure is an outer Nijenhuis tensor $N$ for
a Dirac subbundle $L$ such that "$\s_N$" is skew-symmetric on $Sec(L)$, so
that "$\s_N$" is a deformed Lie algebroid bracket on $L$. A particular
case is when $N$ is a weak Nijenhuis and paired tensor on a Courant
algebroid which is an outer Nijenhuis tensor for $L$. In this case the
subbundle $L$ is a Dirac structure for the deformed Courant algebroid
product "$\s_N$".

Finally, Lie bialgebroids are known as complementary to each other Dirac
subbundles (structures) $E_1,E_2$ in a Courant algebroid $A$, $E_1\op
E_2=A$. It is therefore completely natural to call by \textit{Lie
bialgebroid-Nijenhuis structure} any tensor $N$ on $A$ which yields
Dirac-Nijenhuis structures for both: $E_1$ and $E_2$. The deformed bracket
restricted to $E_1$ and $E_2$ gives two Lie algebroid brackets and the
consistency condition ($N$ is paired) is satisfied, so we get a new Lie
bialgebroid. It is interesting that, associated with particular
contractions, we recover presymplectic-Nijenhuis and Poisson-Nijenhuis
structures (cf. \cite{MM,KSM}). Since the latter play a prominent role in
the theory of integrable systems, this discovery supports once more the
conviction on the importance of bi- or double-structures like Lie
bialgebras, Manin triples, Lie bialgebroids, Courant algebroids, etc., in
complete integrability. Note that a close relation of Poisson-Nijenhuis
structures with Lie bialgebroids was observed first by
Y.~Kosmann-Schwarzbach \cite{KS} (see also \cite{GU}).

\section{Contractions of Leibniz algebras and the Courant bracket}
The language of \textit{Leibniz algebras} is very useful in
description of
\textit{Lie bialgebroids} in the sense of K.~Mackenzie and P.~Xu
\cite{MX}. In \cite{CGM} it has been developed the theory of
contractions for binary operations of arbitrary type, so that all this
general theory of contractions can be directly applied to Leibniz
products (or brackets) on sections of a vector bundle $A$, in
particular for Courant algebroids and Lie bialgebras.

\medskip\noindent
\textbf{Definition 1.} A \textit{Leibniz product (bracket)} on a vector
space $\A$ is a bilinear operation "$\s$" satisfying the Jacobi identity
\begin{equation}\label{3}
(X\s Y)\s Z=X\s(Y\s Z)-Y\s(X\s Z)
\end{equation}
for all $X,Y,Z\in\A$. The space $\A$ equipped with a Leibniz product
we call a \textit{Leibniz algebra}.

\medskip\noindent
Remark that Leibniz algebras as non-skew-symmetric generalizations of Lie
algebras were first studied by J.-L.~Loday \cite{Lo} (they are called
sometimes \textit{Loday algebras}) and a major part of (co)homology theory
of Lie algebras was generalized to Leibniz algebras. Let now "$\s$" be a
Leibniz product on the space $\A=Sec(A)$ of sections of a vector bundle
$A$ over $M$ which is local, i.e. which is locally defined by a
bidifferential operator, and let $N:A\ra A$ be a $(1,1)$-tensor over $A$.
According to the general scheme in \cite{CGM}, if the Nijenhuis torsion
(\ref{1}) vanishes, the contracted product (\ref{2}) is a Leibniz product
which is \textit{compatible} with the original one, i.e. $X\s_N Y+\zl X\s
Y$ is a Leibniz product for any $\zl\in\R$. However, we can have the same
under much weaker conditions.

\begin{lem} The products "$\s_N$" and "$\s$" are always compatible in the sense that
\begin{equation}\label{comp}
(X\s_N Y)\s Z-X\s_N(Y\s Z)+Y\s_N(X\s Z)+ (X\s Y)\s_N Z-X\s(Y\s_N
Z)+Y\s(X\s_N Z)=0.
\end{equation}
\end{lem}
\bepf Direct computations with the use of the Jacobi identity (\ref{3}) for
"$\s$". \epf

\begin{theo}\label{t0} The contracted product (\ref{2}) is still Leibniz if and only if the
Nijenhuis torsion (\ref{1}) is a 2-cocycle with respect to the Leibniz
cohomology operator, i.e.
\bea\label{coc}
(\zd T_N)(X,Y,Z)&=&T_N(X,Y\s Z)-T_N(X\s Y,Z)-T_N(Y,X\s Z)\\
&&-T_N(X,Y)\s Z+X\s T_N(Y,Z)-Y\s T_N(X,Z)=0.\nn
\eea
In this case "$\s_N$" and "$\s$" are compatible Leibniz products.
\end{theo}
\bepf One proves that
\be\label{dd}(X\s_N Y)\s_N Z-X\s_N(Y\s_N Z)+Y\s_N(X\s_N Z)=(\zd T_N)(X,Y,Z)
\end{equation}
by direct computations using the Jacobi identity for "$\s$" and the
compatibility condition (\ref{comp}). In the case when "$\s_N$" is a
Leibniz product, the Jacobi identity for the product $X\s_N Y+\zl X\s Y$
reduces to (\ref{comp}). \epf

\medskip\noindent
The tensor $N$ we will call a \textit{Nijenhuis tensor} (for the
Leibniz algebra $\A$) if the Nijenhuis torsion $T_N$ vanishes and a
\textit{weak Nijenhuis tensor} if the Nijenhuis torsion $T_N$ is a
Leibniz 2-cocycle. In both cases the contracted product "$\s_N$" is
Leibniz and it is compatible with the original one.

An interesting example of a Leibniz product is the following version of
the \textit{Courant bracket} on sections $X+\zx$ of the bundle $\sT M\op
\sT^*M$:
\begin{equation}\label{Cou}(X+\zx)\s(Y+\zh)=[X,Y]+(\Ll_X\zh-i_Yd\zx).
\end{equation}
This is an example of a \textit{Courant algebroid} associated with the
trivial Lie bialgebroid $((\sT M,[\cdot,\cdot]),(\sT^*M,0))$ with the
standard Lie algebroid structure on $\sT M$ and the trivial one on
$\sT^*M$ (cf. \cite{LWX,Ro}). If we have a Nijenhuis tensor $N_0$ for $\sT
M$, we can contract the standard bracket of vector fields to a Lie
algebroid bracket $[X,Y]_{N_0}=[N_0X,Y]+[X,N_0Y]-N_0[X,Y]$ (cf.
\cite{KSM,CGM}). We obtain another trivial Lie bialgebroid $((\sT
M,[\cdot,\cdot]_{N_0}),(\sT^*M,0))$ with the corresponding Courant bracket
\begin{equation}\label{Cou1}(X+\zx)\s^{N_0}(Y+\zh)=[X,Y]_{N_0}+
(\Ll^{N_0}_X\zh-i_Yd^{N_0}\zx),
\end{equation}
where $d^{N_0}$ and $\Ll^{N_0}$ denote the de Rham differential and
the Lie derivative, respectively, associated with the Lie algebroid
$(\sT M,[\cdot,\cdot]_{N_0})$. It is a matter of standard calculations
to show that $d^{N_0}=i_{N_0}d-di_{N_0}$, where $i_{N_0}$ is the
derivation of the algebra of differential forms generated by $N_0$
(see \cite{KSM,GU}). We may as well speak of the product (\ref{Cou1})
purely formally, not even assuming that $N_0$ is a Nijenhuis tensor,
and get the following
\begin{theo}\label{Ni} The product "$\s^{N_0}$" defined by (\ref{Cou1})
is actually the contracted product "$\s_N$" with
$N(X+\zx)=N_0X-{}^tN_0\zx$, where ${}^tN_0:\sT^*M\ra\sT^*M$ is the
dual map: $\la X,{}^tN_0\zx\ran=\la N_0X,\zx\ran$, i.e.
\bea\label{Cou1a}&(X+\zx)\s^{N_0}(Y+\zh)=\cr &[X,Y]_{N_0}+
(N_0X)\s\zh-X\s({}^tN_0\zh)+{}^tN_0(X\s\zh)-({}^tN_0\zx)\s Y+\zx\s
(N_0Y)+{}^tN_0(\zx\s Y).
\eea
\end{theo}
\bepf We have
\beas\la\Ll^{N_0}_X\zh,Y\ran&=&(N_0X)\la\zh,Y\ran-\la\zh,[X,Y]_{N_0}\ran\cr
&=&(N_0X)\la\zh,Y\ran-\la\zh,[N_0X,Y]+[X,N_0Y]-N_0[X,Y]\ran\cr &=&
\la\Ll_{N_0X}\zh+{}^tN_0(\Ll_X\zh)-\Ll_X({}^tN_0\zh),Y\ran.
\eeas
The rest can be proved analogously
\epf

\medskip\noindent
Since, for $N$ being Nijenhuis, the contracted bracket "$\s^{N_0}=\s_N$"
is clearly a Leibniz bracket, the tensor $N$ is automatically weak
Nijenhuis in this case. On the other hand, what is rather unexpected, the
tensor $N$ is a Nijenhuis tensor for the Courant bracket (\ref{Cou}) only
in very particular and rare cases. Namely, we have the following.

\begin{theo}\label{t2} For a Nijenhuis tensor $N_0:\sT M\ra \sT M$
on a connected manifold $M$, the tensor $N:\sT M\op \sT^*M\ra \sT M\op
\sT^*M$, $N(X+\zx)=N_0X-{}^tN_0\zx$, is a Nijenhuis tensor for the Courant
bracket (\ref{Cou}) if and only if $N_0^2=\zl I$ for certain $\zl\in\R$.
\end{theo}
\bepf Since $T_N$ vanishes on $\sT M$ and on $\sT^*M$ separately, the
vanishing of $T_N$ on $\sT M\op \sT^*M$ is equivalent to the system of
identities
\bea{}^tN_0\Ll_X^{N_0}\zh&=&\Ll_{N_0X}({}^tN_0\zh)\label{h1},\\
{}^tN_0i_Yd^{N_0}\zx&=&i_{N_0Y}d({}^tN_0\zx),\label{h2}
\eea
for all $X,Y\in Sec(\sT M)$ and $\zh,\zx\in Sec(\sT^* M)$. The first
one is equivalent to
$$N_0[N_0X,Y]=[X,N_0Y]_{N_0}$$ for all $X,Y\in Sec(\sT M)$ and, due to vanishing of the
Nijenhuis torsion of $N_0$, to
$$N_0^2[X,Y]=[X,N_0^2Y].$$ Since (\ref{h2}) in the presence of (\ref{h1})
can be replaced by
$$({}^tN_0)^2d\la Y,\zx\ran=d\la Y,({}^tN_0)^2\zx\ran,$$
the proof follows by the following lemma. \epf

\begin{lem} If a $(1,1)$-tensor $K:\sT M\ra\sT M$ on a connected manifold $M$
commutes with the adjoint action of vector fields, i.e.
\begin{equation}\label{a2}K[X,Y]=[X,KY]
\end{equation}
for all $X,Y\in Sec(\sT M)$, then $K=\zl I$ for certain $\zl\in\R$.
\end{lem}
\bepf In local coordinates $(x^i)$ and the corresponding coordinate vector fields $(\pa_i)$
we can write $K(\pa_j)=K^i_j(x)\pa_i$ and, according to (\ref{a2}),
$$[\pa_k,K^i_j(x)\pa_i]=\frac{\pa K^i_j}{\pa x^k}(x)\pa_i=0$$ for all $k,j$
(we use the Einstein's summation convention), so the coefficients
$K^i_j(x)=K^i_j$ are constant. Hence, (\ref{a2}) applied to
$X=x^1\pa_k$, $Y=\pa_1$, gives $K^i_k=\zd^i_kK^1_1$, i.e. $K=\zl I$,
where $\zl=K^1_1$. This locally defined constant $\zl$ serves for the
whole $M$, since $M$ is connected. \epf

\medskip\noindent
Note that $(1,1)$-tensors $N_0:\sT M\ra\sT M$ with $N_0^2=\zl I$ and
constant rank are \textit{special} in the terminology of \cite{BC}.
They are proportional to such tensors with $\zl=0,\pm 1$. The case
$\zl=-1$ is the case of an \textit{almost complex structure}, $\zl=1$
is the case of an \textit{almost product structure}, and $\zl=0$ is
the case of an
\textit{almost tangent structure}. If $N_0$ is additionally a
Nijenhuis tensor, we deal with a \textit{complex, product, and
tangent structure}, respectively, cf. \cite{BC}.

\begin{cor} A Nijenhuis tensor $N_0:\sT M\ra\sT M$ gives rise to a
Nijenhuis tensor $N=N_0\op(-{}^tN_0):\sT M\op \sT^*M\ra \sT M\op
\sT^*M$ for the standard Courant bracket if and only if $N_0$ is
proportional to a complex, a product, or a tangent structure on $M$.
\end{cor}
Such structures are extremely interesting from the geometric point of
view. However, from an algebraic point of view, the contracted Courant
brackets for complex and product structures are isomorphic with the
original Courant bracket. To enrich the family of contracted brackets
we will work also with weaker versions of Nijenhuis tensors. This
approach will be systematically developed in the next sections for the
general Courant algebroids.

\section{Contractions of Courant algebroids. Dirac-Nijenhuis structures}
A Courant algebroid is not only a Courant product "$\s$" on sections of a
vector bundle $A$ but also a nondegenerate symmetric pairing
$\la{\cdot},{\cdot}\ran$ on $A$ with certain consistency relations. The
general contraction procedure in such a case is obvious: we contract the
product and check if the consistency conditions with other structures are
still satisfied. If this is the case, we call such contraction the
contraction of the whole structure and the corresponding Nijenhuis tensor
we call the Nijenhuis tensor for the global structure.

Let us recall briefly the structure of a Courant algebroid. We will
use the Leibniz bracket version of the Courant product (bracket)
presented in \cite{Ro} with some simplifications (cf. \cite[Definition
1]{GM}, \cite[Definition 2.1]{KS1} and \cite{Uch}). Thus the
`compressed' definition is as follows.

\medskip\noindent
\textbf{Definition 2.} A \textit{Courant algebroid} is a vector bundle
$\zt:A\ra M$ with a Leibniz product (bracket) "$\s$" on $Sec(A)$, a vector
bundle map (over the identity) $\zr:A\ra \sT M$ and a nondegenerate
symmetric bilinear form $\la{\cdot},{\cdot}\ran$ on $A$ satisfying the
identities
\bea\label{4}
&\zr(X)\la Y,Y\ran=2\la X,Y\s Y\ran,\\ &\zr(X)\la Y,Y\ran=2\la X\s
Y,Y\ran.\label{5}
\eea
Note that (\ref{4}) is equivalent to
\begin{equation}\label{4a}
\zr(X)\la Y,Z\ran=\la X,Y\s Z+Z\s Y\ran.
\end{equation}
Similarly, (\ref{5}) easily implies the invariance of the pairing
$\la{\cdot},{\cdot}\ran$ with respect to the left multiplication
\begin{equation}\label{6}\zr(X)\la
Y,Z\ran=\la X\s Y,Z\ran+\la Y,X\s Z\ran
\end{equation}
and that $\zr$ is the anchor map for the left multiplication:
\begin{equation}\label{7}
X\s(fY)=fX\s Y+\zr(X)(f)Y.
\end{equation}
Assume now that $N$ is a $(1,1)$-tensor on $A$ and consider the
`contracted' product (\ref{2}). We do not assume that $N$ is Nijenhuis at
the moment. Exactly as in the classical case of a Lie algebroid
contraction \cite[Lemma 2]{CGM}, we have the anchor $\zr_N=\zr\circ N$ for
the contracted multiplication
\begin{equation}\label{8}
X\s_N(fY)=f(X\s_NY)+\zr(NX)(f)Y.
\end{equation}
Now, let us check under what conditions the identities (\ref{4}) and
(\ref{5}) are still satisfied for "$\s_N$". Let $N^*$ be the adjoint
of $N$ with respect to the pairing:
$$\la NX,Y\ran=\la X,N^*Y\ran$$
and let $\zD=N+N^*$. Using the invariance (\ref{5}) we get easily
\beas
\la X\s_NY,Z\ran&=&\la NX\s Y+X\s NY-N(X\s Y),Z\ran\\
&=&\zr(NX)\la Y,Z\ran-\la Y,NX\s Z\ran+\zr(X)\la NY,Z\ran-
\la NY,X\s Z\ran-\la X\s Y,N^*Z\ran\\
&=&\zr(NX)\la Y,Z\ran-\la Y,NX\s Z\ran+
\la Y,N^*(X\s Z)\ran+\la Y, X\s N^*Z\ran,
\eeas
which equals $\zr(NX)\la Y,Z\ran-\la Y,X\s_NZ\ran$ if and only if
$$\la Y,X\s \zD Z-\zD(X\s Z)\ran=0$$
for all $X,Y,Z$, i.e. if and only if $\zD$ commutes with the left
multiplication
\begin{equation}\label{9}
X\s \zD Z-\zD(X\s Z)=0.
\end{equation}
Thus (\ref{9}) is equivalent to the invariance of the pairing with respect
to "$\s_N$":
$$\zr_N(X)\la Y,Z\ran=\la X\s_NY,Z\ran+\la Y,X\s_NZ\ran.
$$
Similarly, checking (\ref{4}) for "$\s_N$", we get
\beas
\la X,Y\s_NY\ran&=&\la X,NY\s Y+Y\s NY-N(Y\s Y)\ran\\
&=&\zr(X)\la Y,NY\ran-\la N^*X,Y\s Y\ran\\ &=&\frac{1}{2}\zr(X)\la
Y,\zD Y\ran-\frac{1}{2}\zr(N^*X)\la Y,Y\ran
\eeas
which equals $\frac{1}{2}\zr(NX)\la Y,Y\ran$ if and only if
$$\zr(X)\la Y,\zD Y\ran=\zr(\zD X)\la Y,Y\ran.$$
The latter can be rewritten in the form
$$\la X,Y\s \zD Y+\zD Y\s Y\ran=2\la\zD X,Y\s Y\ran$$
or
$$Y\s \zD Y+\zD Y\s Y=2\zD(Y\s Y).$$
Using (\ref{9}) we get finally the condition
\begin{equation}\label{10}\zD(Y\s Y)=\zD Y\s Y.
\end{equation}

\begin{theo} If $N:A\ra A$ is a $(1,1)$-tensor on a Courant algebroid, then
the contracted product (\ref{2}) is compatible with the symmetric pairing
$\la{\cdot},{\cdot}\ran$ of the Courant algebroid, in the sense that
(\ref{4}) and (\ref{5}) are satisfied for "$\s_N$" and $\zr_N$, if and
only if
$$X\s (N+N^*)
Y=(N+N^*)(X\s Y)\quad \text{and} \quad (N+N^*)(Y\s Y)=(N+N^*)Y\s Y$$ for
all sections $X,Y$ of $A$.
\end{theo}

Of course, how restrictive the above conditions are, depends on `how
irreducible' is the Courant product. However, there is one case which
works for any Courant algebroid, namely the case $N+N^*=\zl I$,
$\zl\in\R$.

\medskip\noindent
\textbf{Definition 3.} A $(1,1)$-tensor on a Courant algebroid we call
\textit{paired} if $N+N^*=\zl I$ for some $\zl\in\R$. A paired (weak)
Nijenhuis tensor we call \textit{(weak) Courant-Nijenhuis tensor}.

\medskip\noindent
Thus weak Courant-Nijenhuis tensors give rise to contractions, or better
to say -- deformations, of Courant algebroids. Note however, that the
structure of a Courant algebroid is extremely rigid and that there are
very few true Courant-Nijenhuis tensors. First, observe that $N$ is a
Courant-Nijenhuis tensor if and only $N-\frac{\zl}{2}I$ is
Courant-Nijenhuis (cf. \cite[Theorem 8]{CGM}), so we can always reduce to
the case when $N+N^*=0$. We have the following generalization of Theorem
\ref{t2}.
\begin{theo}\label{N} If $N$ is a Courant-Nijenhuis tensor with $N+N^*=0$,
then $N^2$ commutes with the left multiplication:
$$X\s N^2Y=N^2(X\s Y)$$
and $N^2(Y\s Y)=(N^2Y)\s Y$.
\end{theo}
\bepf
Using $N^*=-N$ and the invariance of the pairing, we get
\begin{equation}\label{p1}\la N(X\s_NY),Z\ran=-\la X\s_NY,NZ\ran=
-\zr(NX)\la Y,NZ\ran+\la Y,X\s_NNZ\ran
\end{equation}
and
\begin{equation}\label{p2}\la NX\s NY,Z\ran=\zr(NX)\la NY,Z\ran+\la Y,N(NX\s Z)\ran,
\end{equation}
so $N$ is Nijenhuis implies that the r.h. sides of (\ref{p1}) and
(\ref{p2}) are equal, i.e.
\begin{equation}\label{p3}
X\s_NNZ-N(NX\s Z)=0.
\end{equation}
But the l.h.s of (\ref{p3}) is
$$NX\s NZ-N(X\s_NZ)-N^2(X\s Z)+X\s N^2Z$$
and vanishing of the Nijenhuis torsion implies $N^2(X\s Z)=X\s N^2Z$.
The second identity one proves analogously, see the proof of
(\ref{10}).
\epf

\medskip\noindent
\textbf{Remark.} The above property of $N$ is a strong restriction indeed.
We know already that in the case of the standard Courant bracket this
implies that $N^2$ is proportional to the identity (cf. Theorem \ref{t2}).
One can see this problem as the problem of small intersection of the
properties: being paired and being Nijenhuis. Indeed, exactly as in
\cite{CGM}, any Leibniz-Nijenhuis tensor $N$ gives rise to a whole
hierarchy of compatible Leibniz structures and Leibniz-Nijenhuis tensors
of the form $N^k$ while $N^2$, for a paired $N$, is usually not paired.
Thus the concept of a hierarchy for Courant algebroid should be reworked.
For example, one can consider only odd powers or add an additional `twist'
to all powers of $N$. We will not discuss this problem in this note
working, in principle, with generalized versions of Nijenhuis tensors. For
example, one can weaken the assumption for a paired tensor $N$ to
determine a proper contraction assuming just that the tensor $N$ is
weak-Nijenhuis , i.e. we will admit week Courant-Nijenhuis tensors as
well. For a week Courant-Nijenhuis tensor $N$ on a Courant algebroid $A$,
the product "$\s_N$" defines another Courant algebroid product with
respect to the same pairing and the anchor $\zr_N$, and "$\s_N$" is
compatible with "$\s$", i.e. $N+\zl I$ is a one-parameter family of weak
Courant-Nijenhuis tensors (cf. Theorem \ref{t0}).

\medskip\noindent
Let now $L$ be a \textit{Dirac structure} in the Courant algebroid $A$,
i.e. let $L$ be a subbundle which is maximal isotropic and closed with
respect to the Leibniz product "$\s$".

\medskip\noindent
{\bf Definition 4.} The pair $(L,N)$ we call a \textit{Dirac-Nijenhuis
structure} if $N$ is a $(1,1)$-tensor in $A$ such that the deformed
product "$\s_N$" is closed and skew-symmetric on $L$ and the Nijenhuis
torsion $T_N$ vanishes on $L$.

\begin{theo} Let $L$ be a {Dirac structure} in the Courant algebroid
$(A,\s,\la\cdot,\cdot\ran)$
\begin{description}
\item{(a)} If a paired $(1,1)$-tensor $N$ on $A$ is an outer Nijenhuis tensor
for $L$ then $(L,N)$ is a Dirac-Nijenhuis structure.
\item{(b)} If $(L,N)$ is a Dirac-Nijenhuis structure, then $L$
is a Lie algebroid with respect to the product "$\s_N$" and $N(X\s_N
Y)=NX\s NY$ for $X,Y\in Sec(L)$.
\end{description}
\end{theo}
\bepf (a) Since $N$ is paired, the
consistency conditions (\ref{4}), (\ref{5}) are satisfied for "$\s_N$"
that implies the skew-symmetry of "$\s_N$" on any isotropic subbundle.

\smallskip\noindent
(b) The deformed product "$\s_N$" has the anchor $\zr\circ N$ and, due to
(\ref{dd}) the vanishing of the Nijenhuis torsion on $L$ implies that
"$\s_N$" satisfies the Jacobi identity (\ref{3}) on $L$.

\medskip\noindent
\textbf{Examples.} Our Courant algebroid will be $A=\sT M\op\sT^*M$ with
the standard Courant product (bracket)
$$(X+\zx)\s(Y+\zh)=[X,Y]+(\Ll_X\zh-i_Yd\zx).
$$

\smallskip\noindent
\textbf{1.} Let $L$ be the Dirac subbundle in $A$ associated with
a closed 2-form $\zW$, i.e. section of $L$ are of the form $X+\zW
X$ for $X$ being vector fields on $M$. The fact that $\zW$ is
closed can be expressed in terms of the Courant product "$\s$" by
the identity
\be\label{e1}
d\zW(X,Y,\cdot)=X\s\zW Y+\zW X\s X-\zW[X,Y]=0.
\end{equation}
We will refer to any closed 2-form as to a \textit{presymplectic
structure}. Note however that, strictly speaking, a presymplectic
structure is often understood as a closed 2-form of constant rank. We do
not make any assumption on the rank of $\zW$ in this paper. Let $N_0$ be a
$(1,1)$-tensor on $\sT M$ and let $N(X+\zx)=N_0X$ be an associated
$(1,1)$-tensor on $A$. Let us check under what conditions $(L,N)$ is a
Dirac-Nijenhuis structure. First of all, $L$ should be closed with respect
to the deformed bracket "$\s_N$". Since, as easily seen,
\be\label{e0}(X+\zW X)\s_N(Y+\zW
Y)=[X,Y]_{N_0}+N_0X\s\zW Y+\zW X\s N_0Y,
\end{equation}
this condition is equivalent to
\be\label{e2}
N_0X\s\zW Y+\zW X\s N_0Y-\zW[X,Y]_{N_0}=0
\end{equation}
which can be rewritten in the form
\beas
(N_0X\s\zW Y+\zW N_0X\s Y-\zW[N_0X,Y])&&\cr +(\zW X\s N_0Y+X\s\zW
N_0Y-\zW[X,N_0Y])&&\cr-(\zW N_0X\s Y+X\s\zW N_0Y-\zW
N_0[X,Y])&=&\cr d\zW(N_0X,Y,\cdot)+d\zW(X,N_0Y,\cdot)-d(\zW
N_0)(X,Y,\cdot)&=& -d(\zW N_0)(X,Y,\cdot)=0,
\eeas
where we have denoted
$$d(\zW N_0)(X,Y,\cdot)=\zW N_0X\s Y+X\s\zW N_0Y-\zW N_0[X,Y],
$$
independently on the skew-symmetry of $\zW N_0$. But the condition
$$d(\zW N_0)(X,Y,\cdot)=0$$ implies immediately that $\zW N_0$ is
skew-symmetric, i.e. $\zW N_0={}^tN_0\zW$. Indeed,
$$d(\zW N_0)(X,X,\cdot)=d(\zW N_0(X,X))=0$$
for all vector fields $X$, so $\zW N_0(X,X)=0$ for all vector fields $X$
and $\zW N_0$ is skew-symmetric. Thus, $L$ is closed with respect to
"$\s_N$" if and only if $\zW N_0$ is skew-symmetric and $d(\zW N_0)=0$. In
this case
$$(X+\zW X)\s_N(Y+\zW Y)=[X,Y]_{N_0}+\zW[X,Y]_{N_0}.
$$
Finally, the Nijenhuis torsion of $N$ vanishes on $L$ if and only
if
$$N((X+\zW X)\s_N(Y+\zW Y))=N_0([X,Y]_{N_0})=N(X+\zW X)\s N(Y+\zW
Y)=[N_0X,N_0Y],
$$
i.e. $N_0$ is a classical Nijenhuis tensor. This structure is known as
\textit{presymplectic-Nijenhuis structure} (called in \cite{MM} $\zW
N$-\textit{structure}), so that $(L,N)$ as above is a Dirac-Nijenhuis
structure if and only if $(\zW,N_0)$ is a presymplectic-Nijenhuis
structure.

\smallskip\noindent
\textbf{2.} Let $L$ be as above but take the $(1,1)$-tensor on $A$ of the
triangular form: $N(X+\zx)=\zL\zx$, for some $\zL:\sT^*M\ra\sT M$. The
deformed product on $L$ reads
$$(X+\zW X)\s_N(Y+\zW Y)=[X,Y]_{N_0}+N_0X\s\zW Y+\zW X\s N_0Y,
$$
where $N_0=\zL\zW$, so it exactly like (\ref{e0}). We conclude that
$(L,N)$ is Dirac-Nijenhuis in this case if and only if $\zL\zW$ is a
Niejnhuis tensor, $\zW\zL\zW$ is skew-symmetric, and $d(\zW\zL\zW)=0$. In
\cite{MM} such structures are called $\zL\zW$-\textit{structures}.

\smallskip\noindent
\textbf{3.} Let now the Dirac subbundle $L$ of $A$ will be associated with
a Poisson tensor $\zL$, i.e. sections of $L$ are of the form $\zL\zx+\zx$
for $\zx$ being 1-forms, and the Lie algebroid bracket reads
$$(\zL\zx+\zx)\s(\zL\zh+\zh)=[\zL\zx,\zL\zh]+[\zx,\zh]^\zL,
$$
where
$$[\zx,\zh]^\zL=\zL\zx\s\zh+\zx\s\zL\zh$$
is the well-known bracket of 1-forms associated with the Poisson tensor
$\zL$. Put $N(X+\zx)=N_0X$ for some $(1,1)$-tensor $N_0$ on $\sT M$. Since
$$(\zL\zx+\zx)\s_N(\zL\zh+\zh)=[\zL\zx,\zL\zh]_{N_0}+N_0\zL\zx\s\zh+\zx\s
N_0\zL\zh,
$$
requiring the skew-symmetry of this product, we immediately get that
$N_0\zL$ must be skew-symmetric, i.e.
\be\label{PN0}N_0\zL=\zL{}^tN_0,
\end{equation}
and that
$$(\zL\zx+\zx)\s_N(\zL\zh+\zh)=[\zL\zx,\zL\zh]_{N_0}+[\zx,\zh]^{N_0\zL}.
$$
Using (\ref{PN0}) we can rewrite $[\zL\zx,\zL\zh]_{N_0}$ as
$\zL([\zx,\zh]^\zL_{{}^tN_0})$, where $[\cdot,\cdot]^\zL_{{}^tN_0}$ is the
deformation of $[\cdot,\cdot]^\zL$ by ${}^tN_0$, so that the condition
that "$\s_N$" is closed on $L$ can be written as
\be\label{PN}\zL([\zx,\zh]^\zL_{{}^tN_0}-[\zx,\zh]^{N_0\zL})=0.
\end{equation}
The vanishing of the Nijenhuis torsion of $N$ on $L$ takes the form
\be\label{NP}
[\zL\zx,\zL\zh]_{N_0}=[N_0\zL\zx,N_0\zL\zh].
\end{equation}
This simply means that the Nijenhuis torsion of $N_0$ vanishes on the
image of $\zL$. The conditions (\ref{PN0}), (\ref{PN}), and (\ref{NP})
form a weaker version of what is called a \textit{Poisson-Nijenhuis
structure} ($\zL N_0$-\textit{structure} in the terminology of \cite{MM})
for which the conditions are: $N_0\zL$ is skew-symmetric, $N_0$ is
Nijenhuis and (instead of (\ref{PN}))
$$[\zx,\zh]^\zL_{{}^tN_0}-[\zx,\zh]^{N_0\zL}=0
$$
(cf. \cite{MM,KSM}).

\section{Contractions of Lie bialgebroids}

The origin of the concept of Courant algebroid \cite{LWX} was an attempt
to obtain double objects for Lie bialgebroids in the sense of Mackenzie
and Xu \cite{MX}. Suppose now that both $E$ and $E^*$ are Lie algebroids
over $M$ with brackets $[{\cdot},{\cdot}]_E$ and
$[{\cdot},{\cdot}]_{E^*}$, anchors $a$ and $a_*$, respectively. Let $d_E$
(resp., $d_{E^*}$) be the de Rham differential and $\Ll^E$ (resp.,
$\Ll^{E^*}$) be the corresponding Lie derivative associated with the Lie
algebroid structure on $E$ (resp., $E^*$). We will denote sections of $E$
by capitals and sections of $E^*$ by Greek letters and we will often
suppress the indexes in the brackets, de Rham differentials and Lie
derivatives if it will be clear from the context which Lie algebroid they
come from.

On $A=E\op E^*$ there is a natural symmetric nondegenerate bilinear
form:
\begin{equation}\label{11}
\la X+\zx,Y+\zh\ran=\la\zx,Y\ran+\la\zh,X\ran.
\end{equation}
It is well known (cf. \cite[Example 2.6.7]{Ro}) that the bundle $A$
with the symmetric pairing $\la{\cdot},{\cdot}\ran$, the anchor
$\zr=a+a_*$, and the product
\begin{equation}\label{12}
(X+\zx)\s(Y+\zh)=([X,Y]+\Ll_\zx Y-i_\zh dX)+([\zx,\zh]+\Ll_X\zh-i_Y
d\zx)
\end{equation}
is a Courant algebroid if and only if the pair $(E,E^*)$ is a Lie
bialgebroid. The subbundles $E$ and $E^*$ are in this case
\textit{Dirac subbundles}, i.e. maximal isotropic with respect to
the symmetric pairing and closed with respect to the Courant bracket,
transversal to each other. Conversely (see \cite{LWX}), if $L_1$ and
$L_2$ are Dirac subbundles transversal to each other of a Courant
algebroid $A$, then $(L_1,L_2)$ is a Lie bialgebroid, where the
brackets and anchors are just restrictions of the corresponding
structures of the Courant algebroid and $L_2$ is considered as the
dual bundle of $L_1$ under the Courant pairing. The Courant product is
then of the form (\ref{12}) and it is completely determined by the Lie
algebroid structures on $E$ and $E^*$. We have namely
\bea\label{c1}
\la X\s \zh,Y\ran&=&a(X)\la\zh,Y\ran-\la\zh,X\s Y\ran,\\
\la X\s \zh,\zx\ran&=&-a_*(\zh)\la X,\zx\ran+
a(\zx)\la X,\zh\ran+\la X,\zh\s \zx\ran.\label{c2}
\eea
This nice characterization of Lie bialgebroids allows us to define
naturally a concept of contraction of a Lie bialgebroid.

\medskip\noindent
\textbf{Definition 5.} Let $(E,E^*)$ be a Lie bialgebroid and $N$ be a
paired $(1,1)$-tensor on the Courant algebroid $(A=E\op
E^*,\zr,\s,\la{\cdot},{\cdot}\ran)$. The triple $(E,E^*,N)$ we call
\textit{Lie bialgebroid-Nijenhuis structure} if $N$ is an outer Nijenhuis tensor
for both: $E$ and $E^*$.

\begin{theo} If $(E,E^*,N)$ is a Lie bialgebroid-Nijenhuis structure, then
$((E,(\s_N)_{\mid E}),(E^*,(\s_N)_{\mid E^*}))$ is again a Lie
bialgebroid. Moreover, $N$ is a weak Courant-Nijenhuis tensor in the
Courant algebroid $E\op E^*$ and $\s_N$ coincides with the Courant
product $\s^N$ associated with the contracted Lie bialgebroid
$((E,(\s_N)_{\mid E}),(E^*,(\s_N)_{\mid E^*}))$.
\end{theo}
\bepf The contractions $((E,(\s_N)_{\mid E})$ and $(E^*,(\s_N)_{\mid E^*}))$
are clearly Lie algebroid structures on $E$ and $E^*$ respectively.
The tensor $N$ being paired respects the consistency conditions, so
that $((E,(\s_N)_{\mid E}),(E^*,(\s_N)_{\mid E^*}))$ is a Lie
bialgebroid and $\s_N$ is a new Courant bracket, so $N$ is weak
Courant-Nijenhuis tensor. The product $\s_N$ must coincide with
$\s^N$, since the Courant bracket in $E\op E^*$ is uniquely determined
by the Lie algebroid structures in $E$ and $E^*$.
\epf

\medskip\noindent
Let us look closer at the contractions of Lie bialgebroids. First of
all, the splitting $A=E\op E^*$ induces the matrix form of $N$:
\begin{equation}\label{13}N=
\left(\begin{array}{cc}N_E&\zL\\\zW&N_{E^*}\end{array}\right),
\end{equation}
where $N_E$ and $N_{E^*}$ act on $E$ and $E^*$, respectively, and
$\zL:E^*\ra E$, $\zW:E\ra E^*$. The tensor $N$ being paired satisfies
$N+N^*=\zl I$. For $X,Y\in Sec(E)$ we have
$$\la N_E(X)+\zW(X),Y\ran=\la NX,Y\ran=\la X,\zl Y-N_E(Y)-\zW(Y)\ran,$$
so
\begin{equation}\label{14}
\la\zW(X),Y\ran=-\la X,\zW(Y)\ran,
\end{equation}
i.e. $\zW$ is skew-symmetric and can be understood as a section of
$\bigwedge^2E^*$. We will refer to $\zW$ as to a two-form. Similarly,
$\zL$ is a section of $\bigwedge^2E$, referred to as a bivector field.
Finally, it is easy to see that
\begin{equation}\label{15}N_E+{}^tN_{E^*}=\zl I_E,
\end{equation}
where the tensor ${}^tN_{E^*}$ represents the map ${}^tN_{E^*}:E\ra E$
dual to $N_{E^*}:E^*\ra E^*$. Conversely, if $\zL,\zW$ are skew-symmetric
and $N_E$ and $N_{E^*}$ satisfy (\ref{15}), then (\ref{13}) is a paired
tensor.

Clearly, $X\s_N Y=X\s_{N_E}Y+X\s_\zW Y$. Using the obvious notation
$(X+\zx)_E=X$ and $(X+\zx)_{E^*}=\zx$, we get
$$
X\s_NY=X\s_{N_E}Y+(\zW X\s Y+X\s\zW Y)_{E}+(\zW X\s Y+X\s\zW Y-\zW(X\s
Y))_{E^*}.
$$
Thus the condition that $E$ is closed with respect to $\s_N$ reads
\begin{equation}\label{d7}(\zW X\s Y+X\s\zW Y-\zW(X\s Y))_{E^*}=0.\end{equation} But
$$(\zW X\s Y+X\s\zW Y-\zW(X\s Y))_{E^*}=\Ll_X(\zW Y)-i_Yd(\zW X)
-\zW([X,Y]_E)=d\zW(X,Y,{\cdot}),
$$
so that $E$ is closed with the bracket $\s_N$ if and only if $\zW$
is a closed two-form. The analogous statement is, of course, valid
for $E^*$. Note that we will denote the l.h.s of (\ref{d7}) also
$d\zW$ even in the case when $\zW$ is not skew-symmetric. Of
course, in this case $d\zW$ has a meaning as a map and not as a
3-form. Similarly, let us see that
\begin{equation}\label{d8}(\zW X\s Y+X\s\zW
Y)_{E}=\Ll_{\zW X}(Y)-\Ll_{\zW Y}(X)+d_{E^*}(\zW(X,Y))
=[X,Y]^\zW
\end{equation}
is the standard form of the bracket $[{\cdot},{\cdot}]^\zW$ defined on
$E$ by the `bivector field' $\zW\in Sec(\bigwedge^2E^*)$. In the case
when $\zW$ is a `Poisson tensor', i.e. the Schouten bracket
$[\zW,\zW]_{E^*}$ vanishes, the bracket $[{\cdot},{\cdot}]^\zW$ is
known to be a Lie algebroid bracket. We will denote the r.h.s. of
(\ref{d8}) by $[X,Y]^\zW$ also when $\zW$ is not Poisson and not even
skew-symmetric. We get the following.

\begin{theo} Let $(E,E^*)$ be a Lie bialgebroid and let $N$
be a paired tensor of the form (\ref{13}) on the Courant algebroid $E\op
E^*$. Then the subbundle $E$ (resp., $E^*$) is closed with respect to the
contracted bracket "$\s_N$" if and only if $\zW$ (resp., $\zL$) is a
closed two-form with respect to the Lie algebroid structure on $E$ (resp.,
$E^*$), i.e. $\zW\in Sec(\bigwedge^2E^*)$ and $d_E\zW=0$ (resp., $\zL\in
Sec(\bigwedge^2E)$ and $d_{E^*}\zL=0$). In this case the bracket "$\s_N$"
on $E$ (resp., on $E^*$) is of the form $X\s_N Y=[X,Y]_{N_E}+[X,Y]^\zW$
(resp., $\zh\s_N\zx=[\zh,\zx]_{N_{E^*}}+[\zh,\zx]^\zL$).
\end{theo}
Let us now check what means the vanishing of the Nijenhuis torsion on $E$
(and, by duality, on $E^*$). Comparing the parts in $E$ and $E^*$, we get
two equations
\bea
N_E([X,Y]_{N_E}+[X,Y]^\zW)&=&[N_EX,N_EY]_E+(\zW X\s N_EY+N_EX\s\zW
Y)_E,\\ \zW([X,Y]_{N_E}+[X,Y]^\zW)&=&[\zW X,\zW Y]_{E^*} +(\zW X\s
N_EY+N_EX\s\zW Y)_{E^*}.
\eea
They can be rewritten in the form
$$
T_{N_E}(X,Y)+[X,Y]^\zW_{N_E}-[X,Y]^{\zW N_E}=0,
$$
$$
[\zW X,\zW Y]_{E^*}-\zW([X,Y]^\zW)-d(\zW N_E)(X,Y,{\cdot})=0,
$$
where $T_{N_E}$ is the Nijenhuis torsion of $N_E$ with respect to the Lie
algebroid bracket on $E$, the bracket $[{\cdot},{\cdot}]^\zW_{N_E}$ is the
contraction of $[{\cdot},{\cdot}]^\zW$ with respect to $N_E$, the bracket
$[{\cdot},{\cdot}]^{\zW N_E}$ is given by (\ref{d8}) but for (possibly
non-skew-symmetric) $\zW N_E$, and the exterior derivative $d(\zW N_E)$ is
given by (\ref{d7}) but for (possibly non-skew-symmetric) $\zW N_E$. Thus
we get the following.
\begin{theo} The matrix (\ref{13}) acting on $A=E\op E^*$ gives rise to
a Lie bialgebroid-Nijenhuis structure if and only if the following
conditions are satisfied:
\begin{enumerate}
\item $N_E+{}^tN_{E^*}=\zl I_E$ for some $\zl\in\R$;
\item $\zW$ and $\zL$ are skew-symmetric and closed: $d_E(\zW)=0$,
$d_{E^*}(\zL)=0$;
\item The following identities hold:
\be
T_{N_E}(X,Y)+[X,Y]^\zW_{N_E}-[X,Y]^{\zW N_E}=0;\end{equation}
\be\label{v2} [\zW X,\zW Y]_{E^*}-\zW([X,Y]^\zW)-d_E(\zW N_E)(X,Y,{\cdot})=0;
\end{equation}
\be\label{v3}
T_{N_{E^*}}(\zh,\zx)+[\zh,\zx]^\zL_{N_{E^*}}-[\zh,\zx]^{\zL
N_{E^*}}=0;\end{equation}
\be [\zL \zh,\zL\zx]_E-\zL([\zh,\zx]^\zL)-d_{E^*}(\zL N_{E^*})(\zh,\zx,{\cdot})=0.
\end{equation}
\end{enumerate}
\end{theo}
{\bf Remark.} The tensors $\zW N_E$ and $\zL N_{E^*}$ need not be skew
symmetric in general. However, if the Lie algebroid structure on $E$ is
(locally) non-degenerate in the sense that the anchor map, thus $d_E$, is
(locally) non-zero, then they have to be skew-symmetric. Indeed,
(\ref{v2}) implies that $d_E(\zW N_E)(X,X,{\cdot})=0$. But $d_E(\zW
N_E)(X,X,{\cdot})=d_E(\zW(X,X))$, so $\zW(X,X)=0$ and $\zW$ is (locally)
skew-symmetric. Similarly, (\ref{v3}) implies that $[\zh,\zh]^{\zL
N_{E^*}}=0$. But $[\zh,\zh]^{\zL N_{E^*}}=d_E(\zL N_{E^*}(\zh,\zh))$, so
$\zL N_{E^*}(\zh,\zh)=0$ and $\zL N_{E^*}$ is (locally) skew-symmetric.

\medskip
Now consider the trivial Lie bialgebroid $(E,E^*)=(\sT M,\sT^*M)$ with
the standard bracket of vector fields on $\sT M$ and the trivial
bracket on $\sT^*M$. Then $d_{E^*}=0$, $\Ll^{E^*}=0$, the brackets
generated by $\zW$ and $\zW N_E$ are trivial and the above conditions
for the matrix
\begin{equation}\label{17}N=
\left(\begin{array}{cc}\frac{\zl}{2} I+N_0&\zL\\\zW&
\frac{\zl}{2} I-
{}^tN_0\end{array}\right),
\end{equation}
where $\zW$ is a closed 2-form (a presymplectic structure) and $\zL$ is a
bivector field, reduce to
\be\label{w1}
T_{N_0}=0;\end{equation}
\be\label{w2} d(\zW N_0)=0;
\end{equation}
\be\label{w3}
[\zh,\zx]^\zL_{{}^tN_0}-[\zh,\zx]^{\zL{}^tN_0}=0;\end{equation}
\be\label{w4} [\zL \zh,\zL\zx]_E-\zL([\zh,\zx]^\zL)=0.
\end{equation}
Note that, according to the above Remark, in this case $\zW N_0$ and
$\zL{}^tN_0$ are skew-symmetric automatically. The equation (\ref{w1})
means that $N_0$ is a (standard) Nijenhuis tensor which, together with the
presymplectic form $\zW$, constitutes a presymplectic-Nijenhuis structure
($\zW N$-structure) \cite{MM} according to (\ref{w2}). The identity
(\ref{w4}) means that $\zL$ is a Poisson tensor and (\ref{w3}) is a
compatibility condition with $N_0$ which says that we deal with a
Poisson-Nijenhuis structure (cf. \cite{MM,KSM,GU1}). Thus we get the
following.
\begin{theo} The Lie bialgebroid-Nijenhuis tensors $N:\sT M\op \sT^*M
\ra\sT M\op \sT^*M$ for the standard Courant bracket (\ref{Cou}) for the
trivial Lie bialgebroid $(\sT M,\sT^*M)$ are precisely of the form
\begin{equation}\label{17b}N=
\left(\begin{array}{cc}\frac{\zl}{2} I+N_0&\zL\\\zW&
\frac{\zl}{2} I-
{}^tN_0\end{array}\right),
\end{equation}
where $N_0$ is a Nijenhuis tensor, $(N_0,\zW)$ is a
presymplectic-Nijenhuis structure and $(N_0,\zL)$ is a
Poisson-Nijenhuis structure.
\end{theo}
Remark that for a general trivial Lie bialgebroid
$((E,[{\cdot},{\cdot}]),(E^*,0))$ the contracted Lie bialgebroid
associated with the triangular matrix
\begin{equation}\label{17c}N=
\left(\begin{array}{cc}I&\zL\\0&
\ I\end{array}\right),
\end{equation}
is the \textit{triangular Lie bialgebroid} associated with the `Poisson
tensor' $\zL$ in the standard terminology. Note also that the use of outer
Nijenhuis tensors puts a flavor of interaction with the ambient bundle to
the contracted products. For example, the above triangular tensor deforms
the trivial bracket in $E^*$ into a possibly non-trivial bracket
$[{\cdot},{\cdot}]^\zL$ induced by the Lie algebroid structure in $E$.

\section{Concluding remarks}
We have developed the idea of contractions of Courant algebroids, Dirac
structures and Lie bialgebroids as a procedure of deforming such
structures by means of appropriate Nijenhuis tensors. The standard
Nijenhuis tensor approach turned out to be too restrictive, so we had to
deal with tensor whose Nijenhuis torsion vanishes only on a subbundle in
question. We should stress that this idea is of a conceptual nature rather
that an \textit{ad hoc} choice of definitions. The naturality of our
approach is supported by the fact that we can recover basic examples of
the interplay between the fundamental tensors in pairwise dual bundles,
like Poisson-Nijenhuis structures, presymplectic-Nijenhuis structures,
etc., which have been studied in Mathematics and Physics in the context of
integrability. We hope to find direct applications of our formalism in
bihamiltonian formalism and integrability in forthcoming papers.


\end{document}